\documentstyle{amsppt}
\magnification=1100
\NoBlackBoxes
\hsize=13cm
\def\C{\Bbb C}
\def\P{\Bbb P}

\def\Z{\Bbb Z}

\def\ls{\vskip.25in}

\def\ms{\vskip.1in}
\def\Q{\Bbb Q}
\def\L{\Cal L}
\def\O{\Cal O}
\def\y{\bar{y}}

\def\Sym{\text{Sym}}
\topmatter
\title{BEYOND A CONJECTURE OF CLEMENS}\endtitle
\author{ Ziv Ran} \endauthor
\address{Math Department, University of California, Riverside CA 92521}
\endaddress
\email{ziv\@math.ucr.edu}\endemail \abstract We prove some
 lower bounds on  certain twists of the canonical bundle of a
codimension-2 subvariety of a generic hypersurface in projective
space. In particular we prove that the generic sextic
threefold contains no rational or elliptic curves
and no nondegenerate curves of genus 2.\endabstract

\endtopmatter
The geometry of a desingularization $Y^m$ of an arbitrary
subvariety of a generic hypersurface $X^n$ in an ambient variety
$W$ (e.g. $W = \P^{n+1}$) has received much attention over the
past decade or so. Clemens [CKM] has proved that for $m =1 , n =
3, W = \P^4$ and $X$ of degree $d\geq 7$, $Y$ has genus $g \geq
1$, and conjectured that the same is true for $d=6$ as well;
this conjecture is sharp in the sense that hypersurfaces
of degree $d\leq 5$
do contain rational curves.
What the statements mean in plain terms is that,
in the indicated range, a nontrivial function-field solution
of a generic polynomial equation must have genus $>0.$\par
A proof of Clemens' conjecure was published by Voisin [V] who stated
more generally
that, for $X$ of degree $d$ in $\P^{n+1}, n\geq 3, m\leq n-2$,
$p_g(Y)>0$ if $d\geq 2n+1-m$ and $K_Y$ separates generic points if
$d\geq 2n+2-m$. Voisin's proof contained an error, for which the
author subsequently proposed a correction.
In the case of codimension 1, i.e. $m=n-1$,  Xu
[X] gave essentially sharp geometric genus bounds.  For $X^n$ a
generic complete intersection of type $(d_1,...,d_k)$ in any
smooth polarized $(n+k)$-fold $M$, Ein [E] proved that $p_g(Y)>0$
if $d_1+...+d_k\geq 2n+k-m+1$ and $Y$ is of general type if
$d_1+...+d_k\geq 2n+k-m+2.$ Ein's bounds are generally not sharp,
e.g. they fail to yield Clemens' conjecture for $M=\P^{n+1}$. In
[CLR] the authors gave some refinements and generalizations of the
results of Ein and Xu by a method which seems to yield essentially
sharp bounds in codimension 1 but not necessarily in general.
\par
In this paper we give a result which improves and generalizes
Clemens' conjecture.
It gives a kind of lower bound on the canonical bundle
of an arbitrary codimension-2 subvariety
of a generic hypersurface  in an arbitrary ambient
variety. For the case of the sextic in $\P^4$
it shows that the minimal genus of a curve is at least 2,
and at least 3 if the curve is nondegenerate. We proceed
to state the result.
\par
First, we fix a projective space $\P=\P^{n+1}$ and set $L=\O(1)$.
We denote by $\L_d$ the space of homogeneous polynomials
of degree $d$ on $\P$.
\par A
divisor $D$ on a smooth $m$-dimensional variety $Y$ is said to be
{\it{pseudo-effective}} if either $m=1$
and $D$ is effective or $m>1$ and
$D.H^{m-1}\geq 0$ for every ample divisor
$H$ on $Y$. \par

\proclaim{Theorem 1}
Let $X\in\L_d$ be generic, $d\geq n-1\geq 2$, and $f:Y\to X$ a
desingularization of an irreducible subvariety of dimension $n-2$.
Let $p+1$ denote the dimension of the span of  $f(Y)$.
Then either: \par
$$ h^0(K_Y +(n+3-d)f^*L )\geq \min (p,3);\ \  {\text{or}} \tag i $$
$$ K_Y+(n-d-2)f^*L\ \  {\text{is pseudo-effective;\ \ or}}\tag ii$$
(iii)$\ \ \ d\leq n+2$ and
$f(Y)$ is swept out by a family of rational
curves of degree at most $ \frac{2n-2}{d+2-n}$.

\par
Furthermore, if $n=3, d\geq 6$ 
then
$Y$ has genus at least 2.
\endproclaim\par
\remark{Remarks}1. Note that if $n=3$, $Y$ is a curve so there is no
difference between effective and pseudo-effective. Hence in
this case the conclusion (ii) is stronger than (i). Thus
in this case either (i) or (iii) hold. If moreover
$d\geq 6$ then it is well known that $X$ contains no lines,
hence in this case (i) holds. If $n=3, d=6$ we conclude
that $Y$ must have genus at least 3, except possibly if
$Y$ is a genus-2 curve spanning a hyperplane (we
don't believe this case actually occurs).\par
2. We do not know any case where case (ii) holds and
case (i) does not.\endremark\par

For the proof of Theorem 1 we take as our starting point the technique
of [CLR] which 'almost' gives the result, and gives some handle on the
exceptions. To analyze the exceptions we first study subspaces of the
space of
polynomials having a 'large' intersection with the
set of monomials (sect. 1). The proof is
completed in sect.2 . The case $n>3$ uses
bend-and-break, especially Miyaoka's
theorem on generic semipositivity of the cotangent bundle
of a non-uniruled variety .\par
Xi Chen (pers. comm.) has also proven the nonexistence of elliptic
curves on the generic sextic threefold.
His proof is based on degeneration methods.

\heading{1. Linear sections
of monomial varieties}\endheading
 The purpose of this section is to give
a partial classification of linear subspaces of the space of
polynomials which have an improper intersection with the set of
monomials. We begin with a general remark.
\proclaim{Lemma 1.1}
Let $U\subset \P^N$ be an irreducible nondegenerate subvariety
that is the image of a finite morphism from a smooth variety, let
$A\subset\P^N$ be a linear subspace meeting $U$ properly, and let
$B\subsetneq A$ be a linear subspace. Then $B\cap U \subsetneq
A\cap U$ as schemes.\endproclaim
\demo{proof} Let $\pi :U'\to U$
be a finite morphism from a smooth variety and set
$$Z=A\cap U,
Z'=\pi^*(Z).$$ Then $Z'$ is a complete intersection on $U'$ and its
ideal sheaf $I$ admits a resolution by a Koszul complex
$$...\to\oplus\O_{U'}(-2)\to\oplus\O_{U'}(-1).$$ By Kodaira
vanishing, the induced map $\oplus H^0(\O_{U'})\to H^0(I(1))$ is
surjective, which proves the Lemma.\qed\enddemo Note that Lemma
1.1 implies in particular that in the given situation an
{\it{arbitrary}} hyperplane section of $U$ is scheme-theoretically
nondegenerate (hence if reduced, nondegenerate in the usual
sense).\par
Now fix an integer $n\geq 3$ and let $P_k$ denote the
projectivization of the space of homogeneous polynomials of degree
$k$ on $\P^n$, with coordinates $x_0,...,x_n$. Let
$$M_k\subset P_k$$
denote the set of monomials, which may be identified with the
symmetric power $\Sym^k(\P^{n*})$; indeed the composite map
$(\P^{n*})^k\to P_k$ is just the projection of the Segre embedding
$(\P^{n*})^k\to \P^N, N=(n+1)^k-1,$ corresponding to the space of
symmetric tensors. For any $f\in M_{k-j},$ let $$f.M_j\subset M_k$$
be the corresponding 'slice', isomorphic to $M_j$.
\proclaim{Proposition 1.2} Let $K\subset P_k$ be a linear subspace
of codimension $\leq 3$ having an improper intersection with a
general slice of type $f.M_1, f\in M_{k-1}$. Then $K$ has
codimension 3, meets a general $f.M_1$ in codimension 2 and in
fact contains one of the following types of subspaces:\par (A)
$I_L(k)$ for some line $L\subset \P^n$, or\par (B) $I_P^2(k)$ for
some point $P\in\P^n$.\endproclaim
\demo{proof} It was shown in
[CLR], and is in any event easy to see, that any $K$ of
codimension $\leq 2$ meets a general $f.M_1$ properly. We may
therefore assume that $K$ is of codimension 3 and that we have an
irreducible subvariety
$$Z\subset M_k\cap K$$ meeting a general
$f.M_1$ in codimension 2 (hence $Z$ itself has codimension 2).\par
Now consider the natural map $$\pi:(\P^{n*})^k\to M_k =
(\P^{n*})^k/S_k$$ and set $$\sigma^2=\sum_1^k p_i^*(h^2), \sigma
^{1,1}=\sum_{i<j}p_i^*(h)p_j^*(h) \in
H^4_{\Sym}((\P^{n*})^k)=H^4(\Sym^k(\P^{n*})),$$ where $p_i$ are
Cartesian projections and $h\in H^2(\P^{n*})$ is a hyperplane
class. These clearly generate (at least over $\Q$) the symmetric
classes in $H^4((\P^{n*})^k)$. Moreover, it is easy to see by
evaluating on test cycles
$$T'=\pi^{-1}(T)$$ where $T\subset M_k$ is
of the form
$$\P^2.M_{k-1}\ {\text{ or}}\ \P^1_1.\P^1_2.M_{k-2}$$
 for linear
subspaces $\P^2, \P^1_1,\P^1_2\subset \P^{n*}$, that for any
codimension-2 subvariety
$$S\subset M_k,\ {\text{ with}}\ S'=\pi^{-1}(S),$$ we
have $$[S'] = a\sigma^2+2b\sigma^{1,1}, a,b\in\Z_{\geq 0},$$ and
moreover $\sigma^2,2\sigma^{1,1}$ are indivisible, i.e. cannot be
written as
$$m[S'], m>1, S\subset M_k.$$
 Also, if
 $$\alpha\in
H^2_{\Sym}((\P^{n*})^k)$$
is the pullback of a hyperplane class in
$P_k$, then
$$\alpha^2=\sigma^2+\sigma^{1,1}.$$
\par Now returning to
$K$ as above and letting
$$A\subset P_k$$
be any codimension-2
subspace containing $K$, then
$$W=A\cap M_k$$ has codimension 2 and
contains $Z$ and
$$[W']=\sigma^2+2\sigma ^{1,1}.$$
 It follows that
$$[Z']=a\sigma^2+2b\sigma ^{1,1}, a,b\in\{0,1\}.$$
 If $a=b=1$, then
$Z=W$ as cycles, hence as subschemes (as $W$ is pure), which
contradicts Lemma 1.1. Thus
$$[Z']=\sigma^2\  \text{or} 2\sigma^{1,1}.$$ We
assume the latter as the former is similar but simpler. From the
well-known structure of the cohomology of projective spaces it
follows that there exist hyperplanes in $\P^{n*}$, i.e. points
$P_{1,i},P_{2,i}\in\P^n$, so that $$Z'=\sum_{i\neq
j}p_i^*(P_{1,i})\cap p_j^*(P_{2,j}).$$ By symmetry, clearly
$\{P_{1,i},P_{2,i}\}$ is independent of $i$, say $=\{P_1,P_2\}$.
In terms of momomials this says that
$$Z=I_{P_1}(1).I_{P_2}(1).M_{k-2}.$$ Now if $P_1\neq P_2,$ then
$Z$ clearly spans a codimension-2 subspace of $P_k$, which is not
the case. Thus $P_1=P_2=P,$ say, which means precisely that
$Z=I_P^2(k)\cap M_k$, as desired.\qed

\enddemo
{\it {Remark.}\ } It seems likely that this proposition can
be generalized to arbitrary codimension; we hope to
return to this elsewhere.
\ls

\heading{2.Conclusion}\endheading

We now give the proof of Theorem 1.  Let
$f:Y\to X\subset \P=\P^{n+1}$ be a generic member of a filling family as
above and let $N_f, N_{f,X}$ denote the normal sheaves.  Now set
$$D=(\bigwedge^{3}N_f)^{**}.$$ Thus
$$c_1(N_f)=K_Y+(n+2)f^*L\geq D\tag 1$$
where $A\geq B$ means $A-B$ is effective.
Now consider the exact sequence

$$ 0\to N_{f/X}\to N_f\to f^*(L^d)\to 0.$$

Our genericity assumption on $X$ implies that any infinitesimal
deformation of $X$ in $\L_d$ carries an infinitesimal deformation
of $f$, hence the natural map $\L_d\to H^0(f^*L^d)$ admits a
lifting $$\phi:\L_d\to H^0(N_f),\tag 2$$ whence a sheaf map
$$\psi
:\L_d\otimes \Cal O_Y\to N_f,\tag 3$$
such that $H^0(\psi )=\phi$, as well as pointwise maps
$$\psi_y:\L_d\to N_{f,y}.$$
Our filling hypothesis implies
that $f(Y)$ moves, filling up $\P$, hence that $\psi$ is generically
surjective.\par
Now if $H^0(N_{f/X})\neq 0$ (i.e. $(f,Y)$ is not
infinitesimally rigid in a fixed $X$), then arguments of
[CLR] (basically, factoring out the subsheaf generated
by some sections of $N_{f/X}$) already show that case (i)
of Theorem 1 holds. Therefore assume $(f,Y)$ is infinitesimally
rigid in a fixed $X$. Note that this implies that if
$\ G\in\L_d$ is such that $f^*(G)=0,$ then $\psi (G)$
is in the subsheaf $N_{f/X}$, hence $\psi (G)=0.$
Note also that
$\psi$ is uniquely determined.\par
Now arguing as in [CLR], consider the map
$$\psi' :\L\to N,\ \ \ N=N_f$$
obtained by the restriction of $\psi$ on a generic
monomial slice $h_1...h_{d-1}\L$, and
suppose to begin with that
$\psi'$
is generically surjective.
Note that if $h\in\L$ and $\psi'(h)=0$,
then clearly $f^*(h)=0$. Conversely, if $f^*(h)=0$
then $\psi'(h)=0$ by construction.
Also, $\psi'$ clearly
drops rank on $f^*h_1,...,f^*h_{d-1}$, and therefore
we have an exact sequence
$$0\to Z\to \L_Y\otimes\O_Y\to M\to 0,\tag 5$$
where $M\subseteq N$ is a full-rank torsion-free subsheaf with
$$c_1(M)\leq K_Y+(n+3-d)f^*L,$$ and $\L_Y=f^*(\L)\simeq\C^{p+2}$.
This sequence gives rise to a rational map to a Grassmannian, which
by blowing up we may assume is a morphism
$$\gamma:Y\to G:= G(3,p+2),$$
such that $\bigcup\limits_{y\in Y}\gamma (y)$ spans $\C^{p+2}$
and that $c_1(M)=\gamma^*(\O_G(1))$. Thus clearly $c_1(M)$ is
effective and if $p\geq 2$ then $\gamma$ is nonconstant, so
$h^0(c_1(M))\geq 2.$ It is elementary and well known that
any linear $\P^1$ in $G$ consists of the pencil subspaces contained
in a fixed 4-dimensional subspace and containing a fixed 2-
dimensional one, therefore if $p\geq 3$ then $\gamma (Y)$
cannot be contained in such a pencil, so $h^0(c_1(M))\geq 3.$
Thus the conclusions of Part (i) of the Theorem hold\par
Now suppose moreover that $n=3, d=6.$
It suffices to prove that if $Y$ is planar, i.e. $p=1$,
then $g(Y)\geq 2$ (in fact, $\geq 4$). Let
$$I\subset\L_6\times G(2,\P^4)$$ denote
the (open) set of pairs $(X,B)$ where $X$ is a smooth sextic
hypersurface in $\P^4$
and $B$ is a plane in $\P^4,$ and let $J$ denote the set
of pairs $(D,B)$ where $B$ is a plane in $\P^4$ and $D$
is a sextic curve in $B$, not necessarily reduced or irreducible.
Since a smooth sextic cannot contain a plane, there is a natural
morphism
$$\pi:I\to J,$$
$$(X,B)\mapsto (X\cap B,B).$$
Clearly $\pi$ is a fibre bundle. Now a fundamental fact
of plane geometry  [AC] is that the family of reduced irreducible
plane curves of geometric genus $g$ and degree $d$ is of dimension
$3d+g-1.$ It follows easily from this that the locus $J_0\subset J$
consisting of pairs $(D,B)$ such that $D$ is the target of
a nonconstant map from a curve of genus $\leq 3$ is of codimension
$>6$, hence $I_0:=\pi^{-1}(J_0)\subset I$ is also of codimension $>6$
Since $G(2,\P^4)$ is 6-dimensional,
$I_0$ cannot dominate $\L_6$. This proves our assertion.\par

It remains to analyze the case where the restriction of
$\psi$ on a generic
monomial slice fails to be generically surjective,
which we do via Proposition 1.2, letting $K$ denote
the kernel of $\psi_y$ for $y\in Y$
general. This leads to the 2 possible cases (A), (B) as there.
 \par
We analyze case (B) first. Its conclusion, globalized
over $Y$,  may be restated as
follows. Let $P(d)$ be the first principal-parts sheaf (Jet bundle)
of the sheaf $\O(d)$ on the projective space $\P$ associated to $\L$,
whose fibre at a point $P\in\P$ is $\O(d)/I_P^2(d)$
It is a standard fact that
$$P(d)\simeq \L\otimes\O(d-1).\tag 6 $$ What
(B) says is that $\psi$ factors through a map (automatically
generically surjective)
$$f^*(P(d))\to N.$$
 By (6) it follows that
$$N(-(d-1)f^*(L))$$
is generically generated by global sections,
which easily leads to an estimate better than in Theorem 1,
namely $$h^0(K_Y+(n+5-3d)f^*L)>0.$$\par
It remains to analyze case (A). What it yields globally is that there
exists a torsion-free rank-2 quotient $Q$ of the trivial bundle
$\L\otimes\O_Y$ such that the natural map
$$\L\otimes\O_Y\to
f^*(L)$$ factors through $Q$, and such that $\psi$ factors through a
map
$$\gamma:\Sym^dQ\to N.\tag 7$$ Let $R$ be the kernel of the natural map
$Q\to f^*(L)$, which is a torsion-free rank-1 sheaf on $Y$, and
note the diagram
$$\matrix 0&\to & f^*(\Omega_{\P}(1)) & \to &\L\otimes\O_Y
&\to & f^*(L) & \to & 0 \\
           &    & \downarrow &&\downarrow &&\parallel && \\
          0&\to & R & \to & Q & \to & f^*(L) & \to & 0.\endmatrix\tag
          8
          $$

Now let $M\subset\P$ by the line corresponding to to quotient
$\L\to Q_y$  for $y\in Y$ general. We claim next that $M$ is
tangent to $g(Y)$ at $g(y)$. To this end choose coordinates
$x_0,...,x_{n+1}$
on $\P$ so that $x_2,...,x_{n+1}$ vanish on $M$,
$x_0$ vanishes on $f(y)$ and
$x_1$ is general, and let $F$ be the equation of $X$.
Note that, by the rigidity of $(f,Y)$ on a fixed $X$
it follows that as we move $X$ by a projective automorphism,
$(f,Y)$ move by the same automorphism. Consequently it follows that
if $$G=v(F)=\sum h_i\partial/\partial x_i(F)$$
for some global vector field $v\in H^0(T_{\P})$,
then $\psi (G)$ coincides with the normal field to $f$ corresponding
to $v$. As $\psi$ factors through $\gamma$ as above (7),
the value of $\psi (G)$ at $f(y)$ depends only on
the restriction of $G$ on $M$. Now we have
$$d.F|_M=x_0\partial /\partial x_0 (F) + x_1\partial /\partial x_1 (F).$$
Note that
$\psi (F) =0.$ On the other hand as $x_0\partial /\partial x_0$
is a vector field vanishing at $f(y)$, we have
$\psi (x_0\partial /\partial x_0 (F))(g(y))=0.$ Thus
$\psi( x_1\partial /\partial x_1 (F))(g(y))=0$, so the vector
field $x_1\partial /\partial x_1$ is tangent to $f(Y)$ at $f(y)$,
and since $x_1$ doesn't vanish at $f(y)$,
it follows that $M$ is tangent to $f(Y)$ at $f(y)$ , as claimed.\par
It now follows that the map $f^*(\Omega_{\P}(1))\to R$ constructed
above factors through a nonzero map
$$\Omega_Y\otimes f^*(L)\to R.$$
Now if $Y$ is not uniruled, it follows from Miyaoka'a
generic positivity theorem that $R-f^*L$ is pseudo-effective. Fixing as
above a generic monomial $h_1...h_{d-2}$,  $\gamma$ induces a
 map
$$\Sym ^2(Q)\to N$$
which has full rank by, e.g. the main lemma of [CLR] and
which drops rank on $f^*h_1,...f^*h_{d-2}$, hence
$$D-(d-2)f^*(L)-c_1(\Sym ^2(Q))\ \  {\text {is effective}}.$$
Since $c_1(\Sym^2(Q))=3c_1(Q)= 3(f^*(L)+R)$
and $D\leq K_Y+ (n+2)f^*(L),$
the conclusion of part (ii)
clearly holds.\par
Finally suppose that $Y$ is uniruled and let $C$ be the normalization
of a general member of a filling family
of rational curves on $Y$, such that $\Omega_Y|_C$ is seminegative
and
$$K_Y.C\geq -(n-1).$$ This exists by bend-and-break [CKM].
Thus
$$\Omega_Y|_C\simeq \bigoplus\limits_{1}^{n-2} \O_{\P^1}(a_i),\ \ 
{\text{all}}\ a_i\leq 0,\ \ \sum a_i\geq -(n-1).$$
By seminegativity of $\Omega_Y,$ clearly
$$R.C\geq L.C+K_Y.C.$$ Restricting everything
on $C$ yields
$$3(L+R).C\leq (K_Y+f^*(n+4-d)L).C,$$
hence
$$2K_Y.C\leq (n-d-2)L.C.$$ Thus
$$L.C\leq\frac{2n-2}{d+2-n}.$$

Therefore the conclusion of part (iii) holds. \ \qed
\ls\ls
\heading Appendix\endheading
The purpose of this appendix is to give an alternate and
somewhat shorter proof of the main Lemma  (Lemma 1.2) of [CLR],
which goes as follows.
\proclaim{Lemma } Let $Y\subset\P (\L)$
be an irreducible subvariety spanning a $\P^{p+1}.$  Fix integers
$r,k$ with $0\le r-1\le k\le d$. Let $\Cal L_d\to \C^N$ be a
linear map of vector spaces such that for $y_1,\dots,y_{r-1}$
general points of $Y$, the restriction to $\Cal
L_d(-y_1-\dots-y_{r-1})$ induces a surjection
$$\Psi:\Cal
L_d(-y_1-\cdots-y_{r-1})\to\C^{k+1}.$$
Then for a general choice
of elements $h_{k+1},\dots, h_d\in\Cal L$ and for general subsets
$Y_1,\dots,Y_k \subset Y$ each of cardinality $p$, with $Y_i\ni
y_i$, $i=1,\dots,r-1$,  the restriction of $\Psi$ to the subspace
$\Cal L(-Y_1)\cdots \Cal L(-Y_k)h_{k+1}\cdots h_d$ surjects.
\endproclaim
\demo{proof} By a {\it {good chain}} in $\L$ we mean a (connected) 
chain whose components are straight lines (i.e. pencils) of the 
form $\L(-S)$ where $S$ is a general $p$-tuple in $Y$ and whose 
'vertices' (i.e. singular points) are general in $\P (\L )$. Clearly two 
general elements of $\L$ can be joined by a good chain. 
\def\M{\Cal M}
Let $\M_d\subset\L_d$ be the set of monomials. By a {\it {good 
chain}} in $\M_d$ we mean a chain which is a union of subchains of 
the form $\Cal C_i'=h_1\cdots\Cal C_i\cdots h_d$ where $\Cal C_i$ 
is a good chain in $\L$. It is easy to see that two general 
monomials can be joined by a good chain. \par Next, let us say 
that a monomial $h_{d-e+1}\cdots h_d\in\M_e$ is rel $\bar{y}, 
\bar{y}=\{y_1,...,y_{r-1}\}$ if $h_i\in\L(-y_i), i=d-e+1,...,r-1$ 
(this condition is vacuous if $r-1<d-e+1$); denote by 
$\M_e(-\bar{y})$ the set of these. Again it is easy to see that 
two general elements of $\M_e(-\bar{y})$ can be joined by a good 
chain within  $\M_e(-\bar{y})$.\par Now to prove the Lemma it 
suffices to show by induction on $q, 0\le q\leq k$, that, with the 
above notations, 

$$ \dim\Psi (\L(-Y_1)\cdots\L(-Y_q)h_{q+1}\cdots h_d)\geq\min 
(q+1, k+1)$$

(for general choices rel $\bar{y}$). For $q=0$ this is clear as 
$\Psi$ is nonzero on a general element of $\M_d(-\bar{y})$, 
because these span $\L_d(-\bar{y}$). Assume it is true for $q$ and 
false for $q+1$, and suppose first that $q+1\leq r-1$. Now because 
$$Z:=\Psi (\L(-Y_1)\cdots\L(-Y_q)h_{q+1}\cdots h_d) =\Psi 
(\L(-Y_1)\cdots\L(-Y_q)\L(-Y_{q+1})h_{q+2}\cdots h_d)$$ by 
assumption (i.e $Z$ doesn't move as $h_{q+1}$ varies in 
$\L(-Y_{q+1})$, and because two general elements of $\L(-y_{q+1})$ 
can be joined by a good chain (whose components are of the form 
$\L(-Y_{q+1})$), it follows that $Z$ is independent of 
$h_{q+1}\in\L(-y_{q+1})$, fixing the other $h$'s. Since 
$y_1,...y_{r-1}$ are all interchangeable, it follows that a 
similar statement holds for any permutation of them. In particular 
$Z$ contains all pencils of the form $\Psi 
(h_1\cdots\L(-Y_j)\cdots h_d ), j=1,...,r-1$ and from connectivity 
of $\M_{r-1}(-\y)$ by good chains it follows that in fact $Z =\Psi 
(\L_{r-1}(-\y)h_r\cdots h_d)$. As for the remaining $h_j$'s, say 
$j=r$, we may pick $y_r\in f^*(h_r)_0$ and apply similar reasoning 
to $y_1,...,y_{r-1},y_r$ in place of $y_1,...,y_{r-1}$ . The 
foregoing argument yields that $Z$ is independent of 
$h_r\in\L(-y_r)$, fixing $y_r$ and the other $h$'s (indeed that 
$Z=\Psi (L_r(-\y-y_r)h_{r+1}\cdots h_d)$. Then another similar 
argument  with good chains of $h_j$'s not fixing $y_j$ yields 
easily that  $Z$ is actually independent of $h_j\in\L $. Now since 
we can connect two general elements of $\M_d(-\bar{y})$ by a good 
chain, we conclude that  $Z=\Psi (\L_d (-\bar {y}))$, which is a 
contradiction. 

The case $q+1\geq r$ is similar but simpler: we may conclude 
directly that $Z$ is independent of $h_{q+1}\in\L(-Y_{q+1})$, 
hence of $h_j\in\L(-Y_j) $ for all $r\leq j\leq k$ and use good 
chains as above to deduce a contradiction.\qed 
\enddemo

{\bf {Acknowledgement}}\ This work was begun while I was visiting
Taiwan under the auspices of the National Center for Theoretical Sciences,
J. Yu, Director. I would like to thank the NCTS for the invitation, and
especially Prof. L.C.Wang and his colleagues
at National Dong-Hwa University for their outstanding hospitailty.
This paper is a continuation of an earlier one [CLR] done
with Luca Chiantini and Angelo Lopez, and I would like to thank them
for what I have learned from them. Thanks are also due to C. Voisin
for pointing out an error in an earlier version of this paper.

\vfill\eject
\centerline{\bf References}
\vskip .5cm
\item{[AC]} Arbarello, E., Cornalba, M.:'A few remarks', Ann Sci.\ \'Ec.\
Norm.\ Sup.\ 16 (1983), 467-488.
\ms
\item{[C]}Clemens, H.:'Curves in generic hypersurfaces', Ann.\ Sci.\ \'Ec.\
Norm.\ Sup.\ 19 (1986), 629-636.
\ms
\item{[CKM]}Clemens, H. et al.:'Higher-dimensional complex geometry',
Ast\'erisque 1988.
\ms
\item{[CLR]} Chiantini, L., Lopez, A.F., Ran, Z.: 'Subvarieties of generic
hypersurfaces in any variety' (math.AG/9901083; Math. Proc.Camb.Phil. Soc.,
to appear).
\ms

\item{[E]} Ein, L.:'Subvarieties of generic complete intersections',
Invent.\ Math.\ 94 (1988), 163-169; II, Math.\ Ann.\ 289 (1991), 465-471.
\ms
\item{[Mi]} Miyaoka, Y.  :  `The Chern classes and Kodaira dimension
of a minimal variety'.  Adv. Studies in Pure Math. {\bf 10} (1981), 449-476.
\ms

\item{[V]} Voisin, C.: 'On a conjecture of Clemens on rational
curves on hypersurfaces', J.\ Differ.\ Geom.\ 44 (1996), 200-213.
{\it Erratum, \ ibid. }\  49, 1998, 601-611.
\ms

\item{[X]} Xu, G.: 'Subvarieties of general hypersurfaces in projective
space', J.\ Differ.\ Geom.\ 39 (1994), 139-172; 'Divisors on generic complete
intersections in projective space', Trans.\ AMS.\ 348 (1996), 2725-2736.

\ls\ls

\enddocument